\documentclass[12pt,a4paper]{amsart}
\usepackage{amsmath,amssymb,amscd}
\input xy
\xyoption{all}

\newcommand{\lra}{\longrightarrow}
\newcommand{\hra}{\hookrightarrow}
\newcommand{\ra}{\rightarrow}

\newcommand{\CC}{\mathbb K}

\newcommand{\PP}{\mathbb P}

\newcommand{\cO}{\mathcal{O}}

\theoremstyle{plain}
\newtheorem{theorem}{Theorem}[section]
\newtheorem{lemma}[theorem]{Lemma}
\newtheorem{proposition}[theorem]{Proposition}
\newtheorem{cor}[theorem]{Corollary}

\newtheorem{rem}[theorem]{Remark}

\begin{document}
\title[Genus Zero]{Coherent Systems of Genus 0\\ 
II: Existence results for $k\ge3$}

\author{H. Lange}
\author{P. E. Newstead}

\address{H. Lange\\Mathematisches Institut\\
              Universit\"at Erlangen-N\"urnberg\\
              Bismarckstra\ss e $1\frac{ 1}{2}$\\
              D-$91054$ Erlangen\\
              Germany}
              \email{lange@mi.uni-erlangen.de}
\address{P.E. Newstead\\Department of Mathematical Sciences\\
              University of Liverpool\\
              Peach Street, Liverpool L69 7ZL, UK}
\email{newstead@liv.ac.uk}

\thanks{Both authors are members of the research group VBAC (Vector Bundles on Algebraic Curves). They were supported by the Forschungsschwerpunkt 
``Globale Methoden in der komplexen Analysis'' of the DFG. The second author would like to thank the Mathematisches Institut der Universit\"at 
         Erlangen-N\"urnberg for its hospitality}
\keywords{Vector bundle, coherent system, moduli space}
\subjclass[2000]{Primary: 14H60; Secondary: 14F05, 32L10}

\begin{abstract}
In this paper we continue the investigation of coherent systems of type 
$(n,d,k)$ on the projective line which are stable with respect to some value 
of a parameter $\alpha$. We work mainly with $k<n$ and obtain existence results for 
arbitrary $k$ in certain cases, together with complete results for $k=3$. 
Our methods involve the use of the ``flips'' which occur at critical 
values of the parameter.
\end{abstract}
\maketitle

\section{Introduction}

A {\it coherent system of type} $(n,d,k)$ on a smooth projective curve $C$ 
over an
algebraically closed field is by definition a pair $(E,V)$ with $E$ a vector 
bundle of rank $n$ and degree $d$ over $C$ and $V\subset H^0(E)$ a linear
subspace of dimension $k$. For any real number $\alpha$, the 
{\it $\alpha$-slope} 
of a coherent system $(E,V)$ of type $(n,d,k)$ is defined by
$$
\mu_{\alpha}(E,V):=\frac{d}{n}+\alpha\frac{k}{n}.
$$
A {\it coherent subsystem} of $(E,V)$ is a coherent system $(F,W)$ such that 
$F$ is a subbundle of $E$ and $W\subset V\cap H^0(F)$. A coherent system is 
called {\it $\alpha$-stable} ({\it $\alpha$-semistable}) if
$$
\mu_\alpha(F,W)<\mu_\alpha(E,V)\quad (\mu_\alpha(F,W)\le\mu_\alpha(E,V)) 
$$
for every proper coherent subsystem $(F,W)$ of $(E,V)$. According to general 
theory (see, for example, \cite{bgn}), there exists 
a moduli space of $\alpha$-stable coherent systems of 
type $(n,d,k)$, which we denote by $G(\alpha;n,d,k)$. 

In a previous paper \cite{ln}, we obtained necessary conditions for the 
existence of $\alpha$-stable coherent systems of type $(n,d,k)$ on a curve of 
genus $0$. We showed further that these conditions were also sufficient when
$k=1$, but for $k=2$ a special case $(n,d)=(4,6)$ had to be excluded. In this
paper we show that, when $k<n$, the conditions of \cite{ln} remain sufficient 
for the existence of $\alpha$-stable coherent systems for small positive values 
of $\alpha$ (we write this as $\alpha=0^+$). For arbitrary $\alpha$, this is no 
longer true, but we can prove that, for each fixed value of $k$, there are only 
finitely many pairs $(n,d)$ for which exceptional behaviour occurs. When 
$k=3$, there are indeed exceptional cases where the range 
of $\alpha$ for which $\alpha$-stable
coherent systems exist is strictly smaller than the range shown to be 
necessary in \cite{ln}. We analyse these cases and obtain necessary and sufficient 
conditions for existence. We give also an example with $k=4$ to show that, 
in higher ranks,
further complications arise.

We have two principal methods. The first is a development of an argument 
used in \cite{ln}, whereby the existence problem for 
small positive values of $\alpha$ is 
reduced to a problem in projective geometry, which we solve completely. The 
second method is completely different from those of \cite{ln}, depending on an analysis 
of the ``flips'' introduced in \cite{bgn}. The advantage of this approach is 
that it makes it possible to translate results from one value of $\alpha$ to
another. It also allows us to construct $\alpha$-stable coherent systems 
starting from $\alpha$-stable (or even $\alpha$-semistable) coherent systems
of lower rank.

We now outline the content of the paper including statements of the main results 
(for notations, see section 2 or \cite{ln}).
We begin in section 2 by describing the general set-up and establishing
notation. This is followed in section 3 by a detailed
strategy for the analysis of flips. The case where $d$ is a 
multiple of $n$ is considered in section 4, where we prove

\noindent{\bf Theorem 4.5.} {\em Suppose $0<k<n$. 
Then there exists a $0^+$-stable coherent system of type 
$(n,na,k)$ if and only if
$$
ka\ge n-k+\frac{k^2-1}n.
$$}

In section 5, we introduce the concept of an allowable critical data set and 
carry out our first computations of the numbers $C_{12}$ and $C_{21}$ 
associated with the corresponding flips. We prove in particular the following general result:

\noindent{\bf Theorem 5.8.} {\em Let $k$ be a fixed positive integer. Then there are only finitely many 
allowable critical data sets with $n>k$ for which $C_{12}\le0$ or $C_{21}\le0$.}\\

For $k<n$ we write as in \cite{ln}
$$
d=na-t \quad \mbox{and} \quad ka =l(n-k) + t + m
$$
with $0 \leq t <n$ and  $0 \leq m<n-k$.
We then obtain as a consequence of Theorem 5.8:

\noindent{\bf Corollary 5.9.} {\em Let $k$ be a fixed positive integer. Then, for all but finitely many 
pairs $(n,d)$ with $n>k$, one of the following two possibilities holds:
\begin{itemize}
\item $G(\alpha;n,d,k)=\emptyset$ for all $\alpha${\em;}
\item $G(\alpha;n,d,k)\ne\emptyset$ for all $\alpha$ such that 
$$
\frac tk<\alpha<\frac{ln+t}k.
$$
\end{itemize}}

The inequalities for $\alpha$ in this corollary are precisely the necessary conditions of \cite[Propositions 4.1 and 4.2]{ln}. This result therefore
justifies our assertion that, for each value of $k$, there are only finitely 
many pairs $(n,d)$ which exhibit exceptional behaviour. Improved results can 
be obtained when $t=0$ (i.~e. when $d$ is a multiple of $n$), $t=1$ or $t\ge k-1$
(Corollaries \ref{corlge2} and \ref{cort>k}).

In section 6, we reprove the results of \cite{ln}
for $k=2$ using our new techniques. The following three sections contain our 
results for $k=3$. The main result is

\noindent{\bf Theorem 8.4.} {\em Suppose $n \geq 4$. Then
$G(\alpha;n,d,3) \neq \emptyset$ for some $\alpha > 0$ if and only if $l \geq 1$, 
$d \geq \frac{1}{3}n(n-3) + \frac{8}{3}$ 
and $(n,d) \neq (6,9)$.
Moreover, when these conditions hold, $G(\alpha;n,d,3)\ne\emptyset$ if and only if
$$
\frac{t}{3} < \alpha < \frac{d}{n-3} - \frac{mn}{3(n-3)},
$$
except for the following pairs $(n,d)$, where the range of $\alpha$ is as stated :
$$
\begin{array}{cccccccc}
 for& (4,7):&\frac{3}{5} < \alpha < 7;&&&for & (5,9):& \frac{3}{4} < \alpha < \frac{11}{3};\\
 for& (6,11):& 1 < \alpha < \frac{7}{3};&&&for & (7,13):& \frac{3}{2} < \alpha < \frac{8}{3}.\\
\end{array}
$$}

For completeness, we also discuss the case $n\le3$, obtaining

\noindent{\bf Theorem 9.2.} {\em \ 
\begin{itemize}
\item[(i):] $G(\alpha;1,d,3)\ne\emptyset$ if and only if $d\ge2$ and $\alpha>0$.
\item[(ii):] $G(\alpha;2,d,3)\ne\emptyset$ for some $\alpha$ if and only if $d\ge2$. 
Moreover, if $d\ge2$, $G(\alpha;2,d,3)\ne \emptyset$ for all $\alpha>\frac t3$ except 
in the case $d=3$, when $G(\alpha;2,3,3)\ne\emptyset$ if and only if $\alpha>1$.
\item[(iii):] $G(\alpha;3,d,3)\ne\emptyset$ for some $\alpha$ if and only if $d\ge4$. 
Moreover, if $d\ge4$, $G(\alpha;3,d,3)\ne \emptyset$ for all $\alpha>\frac t3$ except 
in the case $d=5$, when $G(\alpha;3,5,3)\ne\emptyset$ if and only if $\alpha>\frac23$.
\end{itemize}}

Finally, in section 10, we give an example of an allowable critical data set 
with $k=4$ where $C_{12}=0$ (this is the smallest value of $k$ for which such 
a critical data set exists).

We work throughout on the projective line $\PP^1$ defined over an 
algebraically
closed field $\CC$.

\section{The set up}

Let $G(\alpha;n,d,k)$ denote the moduli space of $\alpha$-stable coherent 
systems on $\PP^1$ of type $(n,d,k)$. We recall 
\cite[Theorem 3.2]{ln} that, when it is non-empty, $G(\alpha;n,d,k)$ is 
always irreducible of dimension
\begin{equation}\label{eqn:dim}
\beta(n,d,k):=-n^2+1-k(k-d-n).
\end{equation}
In particular, if $G(\alpha;n,d,k)\ne\emptyset$, then
$\beta(n,d,k)\ge0$.

 In accordance
with \cite[section 6]{bgn}, we consider exact sequences
\begin{equation}\label{eqn:ext1}
0 \ra (E_1,V_1) \ra (E,V) \ra (E_2,V_2) \ra 0
\end{equation} and
\begin{equation}\label{eqn:ext2}
0 \ra (E_2,V_2) \ra (E',V') \ra (E_1,V_1) \ra 0
\end{equation}
with $(E,V)$ and $(E',V')$ of type $(n,d,k)$ and $(E_i,V_i)$ of type 
$(n_i,d_i,k_i)$ for $i=1,2$. We suppose also that 
\begin{equation} \label{eqn1.3}
\frac{d_2}{n_2}> \frac{d_1}{n_1} \quad \mbox{and} \quad \frac{k_1}{n_1}>\frac{k_2}{n_2}
\end{equation} 
and define
\begin{equation}\label{eqn:alpha}
\alpha_c = \frac{d_2n - dn_2}{n_2k-nk_2} =\frac{d_2n_1-d_1n_2}{n_2k_1-n_1k_2}.
\end{equation}
We write $\alpha_c^-$ for a value of $\alpha$ slightly smaller than 
$\alpha_c$
and $\alpha_c^+$ for a value of $\alpha$ slightly larger than $\alpha_c$. We 
suppose always that $(E_1,V_1)$ and $(E_2,V_2)$ are both 
$\alpha_c$-semistable.
Note that
$$
\mu_{\alpha_c}(E_1,V_1)=\mu_{\alpha_c}(E_2,V_2)=\frac{d}{n}+\alpha_c\frac{k}{n},
$$
so $(E,V)$ is strictly $\alpha_c$-semistable.

Given this set-up, we shall refer to $\alpha_c$ as a {\it critical value} and
to 
$$A_c:=(\alpha_c,n_1,d_1,k_1,n_2,d_2,k_2)
$$
as a {\it critical data set}. Note that, given $(n,d,k)$, a critical data set
is determined by giving values to $(n_2,d_2,k_2)$ but not necessarily simply
by the critical value $\alpha_c$. Essentially a critical value $\alpha_c$ is
a value of $\alpha$ at which the $\alpha$-stability condition for a coherent
system $(E,V)$ can change as $\alpha$ passes through the value $\alpha_c$,
while the corresponding critical data sets describe the way in which this 
change takes place. For convenience we write 
$G(\alpha_c^-):=G(\alpha_c^-;n,d,k)$ and $G_{\alpha_c}^-$ for the 
``flip locus'' in $G(\alpha_c^-)$, that is the closed subvariety consisting
of those coherent systems which are $\alpha_c^-$-stable but not
$\alpha_c^+$-stable. Similarly we define $G(\alpha_c^+)$ and 
$G_{\alpha_c}^+$ with $+$ and $-$ interchanged.   

As in \cite{bgn} (and putting $g=0$), we define for any critical data set
\begin{equation} \label{eqn1.1}
C_{12} = -n_1n_2 - d_2n_1 + d_1n_2 + k_1(d_2 + n_2 - k_2) 
\end{equation}
and 
\begin{equation} \label{eqn1.2}
C_{21} = -n_1n_2 + d_2n_1 - d_1n_2 + k_2(d_1 + n_1 - k_1).
\end{equation}
We shall explain the significance of $C_{12}$ and $C_{21}$ more precisely in 
section 3. 

 We shall be mainly concerned with the case $0<k<n$. We then write as in the introduction
$$
d=na-t \quad \mbox{and} \quad ka =l(n-k) + t + m
$$
with $0 \leq t <n$ and $0 \leq m<n-k$. Note that, by \cite[Remark 4.3]{ln}, $l>0$ is a necessary condition 
for $G(\alpha;n,d,k)$ to be non-empty.
From (\ref{eqn1.3}), we have $d_2>n_2\frac{d}{n}= n_2a - \frac{n_2}{n}t$ and we write
\begin{equation} \label{eqn3a}
d_2= n_2a+e
\end{equation}
with an integer $e > -\frac{n_2}{n}t$. Using (\ref{eqn3a}), we can rewrite 
(\ref{eqn:alpha}) as
\begin{equation} \label{eqn1.4}
\alpha_c = \frac{ne+n_2t}{n_2k-nk_2}.
\end{equation}
According to \cite[Propositions 4.1 and 4.2]{ln}, we can suppose
\begin{equation}\label{eqn:alphac}
\frac{t}{k} < \alpha_c < \frac{d}{n-k} - \frac{mn}{k(n-k)} = \frac{ln+t}{k},
\end{equation}
which in terms of $e$ means
\begin{equation} \label{eqn4a}
-\frac{k_2}{k}t < e < - \frac{k_2}{k}t + l\left(n_2- \frac{k_2}{k}n\right).
\end{equation}
Note that $-\frac{k_2}{k}t \ge -\frac{n_2}{n}t$ (with equality if and only if
$t=0$), so this inequality is 
stronger than $e>-\frac{n_2}{n}t$. Now write
\begin{equation} \label{eqn1.5}
e= -\frac{k_2}{k}t + l\left(n_2- \frac{k_2}{k}n\right) - \frac{f}{k}
\end{equation}
with an integer $f \geq 1$. In particular
\begin{equation} \label{eqn6}
f \equiv -k_2(t+ln) \equiv k_2m \bmod k.
\end{equation}\\

\section{The strategy}
In this section we explain our strategy for analysing flips. 
The basic idea (introduced in 
\cite{bgn}) is to estimate 
the numbers $C_{12}$ and $C_{21}$ (see (\ref{eqn1.1}) and (\ref{eqn1.2})) 
for any critical data set and use this information to determine how the 
$\alpha$-stability of a coherent system can change as $\alpha$ passes 
through a critical value. We can also use this approach to construct 
$\alpha$-stable coherent systems for values of $\alpha$ close to this 
critical value.

Let $A_c:=(\alpha_c,n_1,d_1,k_1,n_2,d_2,k_2)$ be a critical data set. We 
consider the exact sequences of the forms (\ref{eqn:ext1}), 
(\ref{eqn:ext2}), where as usual we suppose that $(E_1,V_1)$ and $(E_2,V_2)$ 
are both $\alpha_c$-semistable. Our main object in this section is to show
that, in some important cases, the inequalities for the codimensions 
of the flip loci given in \cite[equations (17) and (18)]{bgn} can be 
replaced by equalities. We begin with a version of \cite[Lemma 3.1]{ln} for 
$\alpha$-semistability.

\begin{lemma}\label{lemmass}
Suppose $k>0$ and $(E,V)$ is $\alpha$-semistable for some $\alpha>0$. Then
$$
E\simeq\bigoplus_{i=1}^n{\mathcal O}(a_i)
$$
with all $a_i\ge0$.
\end{lemma}

\begin{proof}
We can write $E=F\oplus G$, where every direct factor of $F$ has negative 
degree and every direct factor of $G$ has non-negative degree. Since $H^0(F)=0$, it follows that
$$
(E,V)=(F,0)\oplus(G,V).
$$
If $F\ne0$, then $\mu_\alpha(F,0)<0$ for all $\alpha$, while 
$\mu_\alpha(G,V)>0$ for $\alpha>0$. This contradicts the 
$\alpha$-semistability of $(E,V)$.
\end{proof}

\begin{cor}\label{cor:ext2=0}
Suppose that $(E_1,V_1)$ and $(E_2,V_2)$ are both $\alpha$-semistable. 
Then
\begin{equation}\label{eqn:ext2=0}
\mbox{\em Ext}^2((E_1,V_1),(E_2,V_2))=\mbox{\em Ext}^2((E_2,V_2),(E_1,V_1))=0.
\end{equation}
\end{cor}
\begin{proof}
This follows from the lemma, together with the formula for 
$\mbox{Ext}^2$ given in 
\cite[equation (11)]{bgn} and the fact that the canonical bundle has 
negative degree.
\end{proof}
\begin{cor}\label{cor:birat} 
Let $\alpha_c$ be a critical value. Suppose that, for every critical data set 
$A_c$ with critical value $\alpha_c$, we have $C_{12}>0$ and $C_{21}>0$. Then 
$G(\alpha_c^+)$ is birational to $G(\alpha_c^-)$. In fact, if $C_{12}>0$ for 
every $A_c$ and $G(\alpha_c^-)\ne\emptyset$, then $G_{\alpha_c}^-$ has positive 
codimension in $G(\alpha_c^-)$. Similarly, if $C_{21}>0$ for 
every $A_c$ and $G(\alpha_c^+)\ne\emptyset$, then $G_{\alpha_c}^+$ has positive 
codimension in $G(\alpha_c^+)$.
\end{cor}
\begin{proof}
Since all non-empty moduli spaces have the expected dimensions (given by 
(\ref{eqn:dim})), it follows from Corollary \ref{cor:ext2=0} and 
\cite[equations (17) and (18)]{bgn} that the flip loci have positive 
codimensions. The result follows.
\end{proof}

The key fact about the numbers $C_{12}$ and $C_{21}$ is that they play two 
r\^oles, as estimates for codimensions of flip loci and for dimensions of 
spaces of extensions. In fact, if we assume in addition to (\ref{eqn:ext2=0})
that

\begin{equation}\label{eqn:hom=0}
\mbox{Hom}((E_1,V_1),(E_2,V_2))=\mbox{Hom}((E_2,V_2),(E_1,V_1))=0,
\end{equation}
then we deduce at once from \cite[equation (8)]{bgn} that
\begin{equation}\label{eqn:c12}
C_{12}=\dim\mbox{Ext}^1((E_1,V_1),(E_2,V_2))
\end{equation}
and
\begin{equation}\label{eqn:c21}
C_{21}=\dim\mbox{Ext}^1((E_2,V_2),(E_1,V_1)).
\end{equation}
In particular, if (\ref{eqn:hom=0}) holds, we always have
$$
C_{12}\ge0,\quad C_{21}\ge0.
$$
 
\begin{lemma}\label{lemmacodim}
Suppose that, for some critical data set $A_c$, there exist 
$\alpha_c$-stable coherent
systems $(E_1,V_1)$ and $(E_2,V_2)$, and that $A_c$ is the only critical
data set for the critical value $\alpha_c$ . Then
\begin{itemize}
\item[(a)] if $C_{21}>0$, the flip locus  $G_{\alpha_c}^-$ is irreducible 
and has codimension  $C_{12}$ in $G(\alpha_c^-)$;
\item[(b)] if $C_{12}>0$, the flip locus  $G_{\alpha_c}^+$ is irreducible 
and has codimension  $C_{21}$ in $G(\alpha_c^+)$.
\end{itemize}\end{lemma} 

\begin{proof}
(a): Consider first the non-trivial extensions (\ref{eqn:ext1}) with  
$(E_1,V_1)$ and $(E_2,V_2)$
both $\alpha_c$-stable. It is easy to see that $(E,V)$ has 
(\ref{eqn:ext1}) as its unique Jordan-H\"older
filtration at $\alpha_c$. Since $\frac{k_1}{n_1}>\frac{k_2}{n_2}$, it follows
also that $(E,V)$ is $\alpha_c^-$-stable. Since $(E_1,V_1)$ and $(E_2,V_2)$ are
non-isomorphic and $\alpha_c$-stable of the same $\alpha_c$-slope, 
(\ref{eqn:hom=0}) holds and therefore also (\ref{eqn:c21}). These 
extensions therefore define 
a non-empty open subset $U$ of $G_{\alpha_c}^-$ of dimension 
$$
\dim U=\dim G(\alpha_c;n_1,d_1,k_1)+
\dim G(\alpha_c;n_2,d_2,k_2)+C_{21}-1.
$$
It follows 
from \cite[Corollary 3.7]{bgn} that $U$ has codimension $C_{12}$ in 
$G(\alpha_c^-)$. 

It remains to show that $G_{\alpha_c}^-$ is irreducible. For this, note  
that, by \cite[Lemma 6.5(ii)]{bgn}, all elements of
$G_{\alpha_c}^-$ come from extensions (\ref{eqn:ext1}) with $(E_1,V_1)$ and 
$(E_2,V_2)$ $\alpha_c^-$-stable. Since $\frac{k_1}{n_1}>\frac{k_2}{n_2}$, 
$\mu_{\alpha_c^-}(E_1,V_1)<\mu_{\alpha_c^-}(E_2,V_2)$. Hence
$\mbox{Hom}((E_2,V_2),(E_1,V_1))=0$, which implies 
(\ref{eqn:c21}). The irreducibility of $G_{\alpha_c}^-$ now follows from 
the irreducibility of the moduli spaces $G(\alpha_c^-;n_1,d_1,k_1)$ and
$G(\alpha_c^-;n_2,d_2,k_2)$.

The proof of (b) is similar.
\end{proof}

\begin{cor}\label{corflip}
Suppose that the hypotheses of the lemma hold. Then one of the following situations 
occurs: 
\begin{itemize}
\item $C_{12}>0$ and $C_{21}>0${\em :} 
$G(\alpha_c^-)$ and
$G(\alpha_c^+)$ are both non-empty and birational to each other;
\item $C_{12}=C_{21}=0${\em :} the flip loci are empty and 
$G(\alpha_c^-)=G(\alpha_c^+)$;
\item $C_{21}=0$, $C_{12}>0${\em :} $G(\alpha_c^-)=\emptyset$,
$G(\alpha_c^+)=G_{\alpha_c}^+\ne\emptyset$; 
\item $C_{12}=0$, $C_{21}>0${\em :} $G(\alpha_c^+)=\emptyset$,
$G(\alpha_c^-)=G_{\alpha_c}^-\ne\emptyset$.
\end{itemize} 
\end{cor}

\begin{proof}
For the first part, the non-emptiness follows from the lemma and the 
birationality is a special case of Corollary \ref{cor:birat}. 
If $C_{12}=C_{21}=0$, then (\ref{eqn:c12}) and (\ref{eqn:c21}) imply that 
the flip loci are empty; this proves the second part. If 
$C_{21}=0$ and $C_{12}>0$, the lemma implies that 
$G(\alpha_c^+)=G_{\alpha_c}^+\ne\emptyset$. It now follows from 
\cite[Corollary 3.4]{ln} that $G(\alpha_c^-)=\emptyset$. The last part is 
proved similarly.
\end{proof}  

\begin{rem}\label{rem1}
{\em In a calculation it may happen that $C_{12}$ or $C_{21}$ comes out to be 
negative. In this case either $(E_1,V_1)$ or $(E_2,V_2)$ fails to exist and
the flip loci are empty.} 
\end{rem}
\begin{rem}\label{rem2}
{\em Suppose there is more than one critical data set $A_c$ for a 
critical value 
$\alpha_c$, such that, for each $A_c$, there exist $\alpha_c$-stable 
coherent systems $(E_1,V_1)$ and $(E_2,V_2)$. We then ignore all $A_c$ for 
which $C_{12}=C_{21}=0$ and replace the remaining $C_{12}$ and 
$C_{21}$ by their minimum
values taken over the various $A_c$. The conclusions of Corollary
\ref{corflip} then hold except possibly when both $C_{12}$ and $C_{21}$ 
have minimum value 
zero (necessarily for different $A_c$). It then follows from the proof of 
Lemma \ref{lemmacodim} that both $G(\alpha_c^-)$ and $G(\alpha_c^+)$ are 
non-empty, but $G(\alpha_c)$ is empty. This contradicts
\cite[Corollary 3.4]{ln}, 
so this situation can never arise.}
\end{rem}
\begin{rem}\label{rem3}
{\em The conclusions of Remark \ref{rem2} still hold if there are additional 
critical data sets with critical value $\alpha_c$, provided these all have 
$C_{12}>0$ and $C_{21}>0$, 
}
\end{rem}

On occasion, we shall need to use extensions in which $(E_1,V_1)$ and 
$(E_2,V_2)$ are not $\alpha_c$-stable. In this paper, it will 
be sufficient to consider extensions
\begin{equation}\label{eqn:nots}
0\ra(\cO(b)^r,0)\ra(E,V) \stackrel{p}{\ra}(E_1,W)\ra0,
\end{equation}
for certain $(E_1,W)$ with $\mu_{\alpha_c}(E_1,W)=b$ (note that 
$\mu_\alpha(\cO(b)^r,0)=b$ for all $\alpha$).

\begin{lemma}\label{lemext}
Suppose that, in {\em(\ref{eqn:nots})}, either $(E_1,W)$ is $\alpha_c$-stable and
\begin{equation}\label{eqn:hyp1}
\dim\mbox{\em Ext}^1((E_1,W),(\cO(b),0))\ge r
\end{equation}
or 
$$
(E_1,W)=\bigoplus_{i=1}^{n-r}(\cO(b'),W_i),
$$
where the $W_i$ are distinct subspaces of dimension $1$ of 
$H^0(\cO(b'))$ and
\begin{equation}\label{eqn:hyp2}
(n-r)b'>r,\quad n-r\ge2.
\end{equation}
Then, for the general extension {\em(\ref{eqn:nots})}, $(E,V)$ is 
$\alpha_c^+$-stable.
\end{lemma}

\begin{proof}
Suppose $(F,U)$ is a subsystem of $(E,V)$ which contradicts 
$\alpha_c^+$-stability. Then $(F,U)$ also contradicts $\alpha_c$-stability.
Since $(E,V)$ is $\alpha_c$-semistable, $(F,U)$ is also 
$\alpha_c$-semistable with the same $\alpha_c$-slope.

In the first case, this implies that the image $p(F,U)$ is either $0$ or 
equal to $(E_1,W)$. If $p(F,U)=0$, then $(F,U)\subset(\cO(b)^r,0)$ and does 
not contradict $\alpha_c^+$-stability of $(E,V)$. So $p(F,U)=(E_1,W)$ and 
we have a diagram
 $$
\xymatrix{
  0 \ar[r] & (\cO(b)^s,0) \ar[r] \ar@{^{(}->}[d] & (F,U) \ar[r] \ar@{^{(}->}[d] & (E_1,W) \ar[r] \ar@{=}[d] & 0\\
  0 \ar[r] & (\cO(b)^{r},0) \ar[r] & (E,V) \ar[r] &  (E_1,W)  \ar[r] & 0       } 
$$ 
with $s < r$. The extensions (\ref{eqn:nots}) are classified by $r$-tuples
$(e_1,\ldots,e_r)$ with $e_i\in\mbox{Ext}^1((E_1,W),(\cO(b),0))$. 
By (\ref{eqn:hyp1}), the general extension (\ref{eqn:nots}) has $e_1,\ldots,e_r$ 
linearly independent. Thus the diagram above is impossible.

In the second case, note first that
$$
\mbox{Hom}\;((\cO(b'),W_i),(\cO(b),0))=0.
$$ 
Hence, by (\ref{eqn:c12}) and (\ref{eqn1.1}),
\begin{equation}\label{eqn:bb}
\dim\;\mbox{Ext}^1((\cO(b'),W_i),(\cO(b),0))=-1-b+b'+(b+1)=b'.
\end{equation}
If $p(F,U)$ is either $0$ or $(E_1,W)$, we argue as in the first case. 
Otherwise note that, since $(\cO(b'),W_i)\not\simeq(\cO(b'),W_j)$ 
for $i\ne j$, there are only finitely many possible choices for $p(F,U)$ and we 
can suppose without loss of generality that  
$$
p(F,U) = \bigoplus_{i=1}^j (\cO(b'), W_i)
$$ 
for some $j$ with $1\le j\le n-r-1$.
If $p$ maps $(F,U)$ isomorphically to $\bigoplus_{i=1}^j (\cO(b'), W_i)$, 
then the extension (\ref{eqn:nots})
restricted to $\bigoplus_{i=1}^j (\cO(b'), W_i)$ splits. But in general this 
does not happen, since by (\ref{eqn:bb}) 
$$
\dim \; \mbox{Ext}^1\left(\bigoplus_{i=1}^j (\cO(b'), W_i), 
(\cO(b),0)\right) = jb' > 0.
$$

We can therefore suppose that the kernel of $(F,U) \ra \bigoplus_{i=1}^j (\cO(b'), W_i)$ has
the form $(\cO(b)^s,0)$ with
$1 \leq s \leq r$. But then $(F,U)$ contradicts $\alpha_c^+$-stability of 
$(E,V)$ 
if and only if
\begin{eqnarray*}
\mu_{\alpha_c^+}(F,U) &= &\frac{bs+b'j}{s+j} + \alpha_c^+ \frac{j}{s+j} 
\\&\ge& \mu_{\alpha_c^+}(E,V)= \frac{br+b'(n-r)}{n} + \alpha_c^+ \frac{n-r}{n}
\end{eqnarray*}
This is equivalent to
$$
(bs +b'j)n -(br+b'(n-r))(s+j) \ge \alpha_c^+ ((n-r)(s+j)-jn),
$$
which in turn is equivalent to
$$(b-b')(s(n-r)-jr)\ge\alpha_c^+(s(n-r)-jr).
$$
Since $b-b'=\alpha_c$, this reduces to 
\begin{equation}\label{eqn:red}
jr \ge s(n-r).
\end{equation}

Now consider the diagram
 $$
\xymatrix{
  0 \ar[r] & (\cO(b)^s,0) \ar[r] \ar@{^{(}->}[d] & (F,U) \ar[r] \ar@{^{(}->}[d] & \bigoplus_{i=1}^j (\cO(b'),W_i) \ar[r] \ar@{=}[d] & 0\\
  0 \ar[r] & (\cO(b)^r,0) \ar[r] \ar@{=}[d] & (F',U') \ar[r] \ar@{^{(}->}[d] & \bigoplus_{i=1}^j (\cO(b'),W_i) \ar[r] \ar@{^{(}->}[d] & 0\\
  0 \ar[r] & (\cO(b)^r,0) \ar[r] & (E,V) \ar[r] &  \bigoplus_{i=1}^{n-r} (\cO(b'),W_i)  \ar[r] & 0       } 
$$ 
where the lower half is the pull-back diagram which always exists, and the 
upper half is a push-out diagram the existence of 
which we have to analyse. The extensions of the middle row are classified by 
$r$-tuples $(e_1,...,e_r)$
with $e_l \in \mbox{Ext}^1\left(\bigoplus_{i=1}^j (\cO(b'),W_i),
(\cO(b),0)\right)$,
which, by (\ref{eqn:bb}), is of 
dimension $jb'$. Hence, for a general extension (\ref{eqn:nots}), such a 
diagram cannot exist unless $jb'\le s$. Combining this with (\ref{eqn:red}),
we obtain
$$
jb'(n-r)\le s(n-r)\le jr,
$$
which contradicts the hypothesis $(n-r)b'>r$. This completes the proof.
\end{proof}

\begin{rem}\label{rem4}
{\em The hypotheses (\ref{eqn:hyp1}) and (\ref{eqn:hyp2}) in 
the statement of the lemma are sharp. In the first case, if 
(\ref{eqn:hyp1}) fails, $(E,V)$ has a direct factor of the form $(\cO(b),0)$ and is not even 
$\alpha_c^+$-semistable. In the second case, if $(n-r)b'\le r$, we can 
take $s=jb'$ to contradict $\alpha_c^+$-stability. In fact, if 
$(n-r)b'=r$, the general $(E,V)$ is strictly $\alpha_c^+$-semistable.}
\end{rem}

\section{The case $t=0$}

In this section we assume $0<k<n$ and consider the existence problem for $0^+$-stable coherent 
systems of type $(n,d,k)$. Note that, if $(E,V)$ 
is such a coherent system, the bundle $E$ is semistable, so 
$E\simeq\cO(a)^n$ and $t=0$. We therefore suppose that $E=\cO(a)^n$ and 
assume also that the homomorphism $\beta:V\otimes\cO\ra\cO(a)^n$ is 
injective. For $1\le q\le k$, we then define
$$
\delta_q(n,a,\beta) = \left\{ \begin{array}{l} \mbox{minimal rank of a direct factor of} \,\, \cO(a)^n \\ 
\mbox{containing the image of some} \,\, \cO^q \subset V\otimes\cO \,\, \mbox{under}\\ \mbox{the composed map} \,\,
\cO^q \hra V\otimes\cO \stackrel{\beta}{\ra} \cO(a)^n. \end{array} \right.
$$

\begin{lemma}\label{lemma:0+}
$(\cO(a)^n,V)$ is $0^+$-stable if and only if $\delta_k(n,a,\beta)=n$ and
$$
\delta_q(n,a,\beta)>\frac{qn}k\quad\mbox{for}\quad1\le q\le k-1.
$$
\end{lemma}

\begin{proof}
Suppose $(F,W)$ is a coherent subsystem of $(\cO(a)^n,V)$ with $\dim W=q$. 
Then $(F,W)$ contradicts the $0^+$-stability of $(\cO(a)^n,V)$ if and only 
if $F\simeq\cO(a)^r$ where either $q=k$ and $r<n$ or $1\le q<k$ and 
$\frac{q}r\ge\frac{k}n$, i.~e. $\frac{qn}k\ge r$. The result follows 
from the definition of $\delta_q$.
\end{proof}

We now convert this condition into a statement in projective geometry. 
For this, let $q,k$ and $n$ denote positive integers with $q \leq k < n$ 
and consider the Segre embedding
$$
\PP^{k-1} \times \PP^{n-1} \hra \PP^{kn-1}.
$$
For any integer $a$ with $0\le a\le kn-2$, let $R(n,a,k,q)$ denote the maximum 
number $r$ 
such that any linear subspace $W \subset \PP^{kn-1}$ of codimension $a+1$ 
contains some
subspace $\PP^{q-1} \times \PP^{r-1} \subset \PP^{k-1} \times \PP^{n-1} \subset \PP^{kn-1}$. 
If $a\ge kn-1$, we define 
$R(n,a,k,q) = 0$. Note that the condition on $W$ is equivalent to saying that $W$ contains the subspace 
$\PP^{qr-1} \subset \PP^{kn-1}$ which is spanned by the image of 
$\PP^{q-1} \times \PP^{r-1}$ in $\PP^{kn-1}$.

\begin{lemma} \label{lem5.1} For a general choice of 
$V\subset H^0(\cO(a)^n)$,
$$
\delta_q(n,a,\beta)=n-R(n,a,k,q).
$$
\end{lemma}

\begin{proof} The map $\beta$ is given by a matrix of the form 
$$
M = \left(f_{ij} \right)_{1\le i\le n,1\le j\le k}          
$$
where the $f_{ij}$ are binary forms of degree $a$. The composition 
$\cO^q \ra \cO(a)^n$ is given
by a matrix $MN_q$ of rank $q$ with 
$$
N_q = \left(b_{jp} \right)_{1\le j\le k,1\le p\le q}
$$
where the $b_{jp}$ are constants and $\mbox{rk}\;N_q=q$.     
By definition of $\delta_q(n,a,\beta)$ we have
$$
n-\delta_q(n,a,\beta) = \max_{A \in GL(n,\mathbb C), \;\; rk \; N_q = q} \,\, \{\mbox{number of zero rows in} \,\, AMN_q \}.
$$
But this equals the maximum number of linearly independent vectors 
$(\lambda_1, \ldots , \lambda_n) \in \mathbb C^n$ such that
$$
(\lambda_1, \cdots, \lambda_n)MN_q = 0,   
$$
the maximum to be taken over all $k \times q$-matrices $N_q$ 
of rank $q$. 

Now let $W$ denote the projectivisation of the kernel of the linear 
map $\mathbb C^{kn} \lra H^0(\cO(a))$ given by 
$$
(\mu_{11},\ldots,\mu_{1k},\ldots,\mu_{n1},\ldots,\mu_{nk}) \mapsto 
\sum f_{ij}\mu_{ij}.
$$
Note that, if $a\le kn-2$, then, for a general choice of the $f_{ij}$, 
$W$ has codimension $a+1$ in $\PP^{kn-1}$. The result follows easily from
the definitions of $\delta_q$ and $R$.
\end{proof} 

The next step is to estimate $R(n,a,k,q)$. 

\begin{lemma} \label{lempg}
$$
R(n,a,k,q) \leq \left\lfloor\frac{1}{2}\left( n-q(a+1) + \sqrt{(n-q(a+1))^2 + 4q(k-q)}\right)\right\rfloor.
$$
\end{lemma}

\begin{proof}
For $a\ge kn-1$, this is obvious since $R(n,a,k,q)=0$. Otherwise,
let $Gr: =Gr(kn-a-1,kn)$ denote the Grassmannian of subspaces of 
codimension $a+1$ in $\PP^{kn-1}$. For a fixed linear subspace 
$\PP^{qr-1} \subset \PP^{kn-1}$, let $\Sigma$ denote the closed subspace 
of $Gr$ consisting of all $W \in Gr$ with $\PP^{qr-1} \subset W$. Finally 
write 
$\Psi: = Gr(q,k) \times Gr(r,n)$. 

We can clearly ignore the values of $r$ for which $\Sigma=\emptyset$, 
or equivalently $a\ge kn-qr$. Otherwise, a necessary condition for a general subspace 
$W \subset \PP^{kn-1}$ of codimension $a+1$ to contain some 
subspace $\PP^{q-1} \times \PP^{r-1} \subset \PP^{k-1} \times \PP^{n-1} \subset \PP^{kn-1}$ is 
$$
\dim \Sigma + \dim \Psi \geq \dim Gr.
$$
Since
$$
\begin{array}{c}
\dim Gr = (kn-a-1)(a+1)\\
\dim \Sigma = (kn-a-1 -qr)(a+1)\\
\dim \Psi = q(k-q) + r(n-r),
\end{array}
$$
this means
$$
(kn-a-1 -qr)(a+1) + q(k-q) + r(n-r) \geq (kn-a-1)(a+1),
$$
which is equivalent to 
$$
r^2 + (q(a+1)-n)r - q(k-q) \leq 0.
$$
This quadratic equation in $r$ always has two real solutions. Solving this 
equation gives the assertion.
\end{proof}

\begin{cor}\label{cor:kk}
$R(n,a,k,k)=0$ if $ka\ge n-k$.
\end{cor}\begin{proof} The hypothesis states that $n-k(a+1)\le0$. The 
assertion then follows immediately from the lemma.
\end{proof}

\begin{theorem}\label{th:0+}
Suppose $0<k<n$. Then there exists a $0^+$-stable coherent system of type 
$(n,na,k)$ if and only if
\begin{equation}\label{eqn:ka}
ka\ge n-k+\frac{k^2-1}n.
\end{equation}
\end{theorem}
\begin{proof}
Note first that (\ref{eqn:ka}) is equivalent to the Brill-Noether inequality 
$\beta(n,na,k)\ge0$ (see (\ref{eqn:dim})). The inequality (\ref{eqn:ka}) is 
therefore a necessary condition for the existence of $\alpha$-stable 
coherent systems of type $(n,na,k)$.

Conversely, suppose (\ref{eqn:ka}) holds. In view of Lemmas \ref{lemma:0+} 
and \ref{lem5.1} and Corollary \ref{cor:kk}, it is sufficient to show 
that, for $1\le q\le k-1$, 
\begin{equation}\label{eqn:R}
R(n,a,k,q)<n-\frac{qn}k.
\end{equation}
We prove first
\begin{lemma}\label{lemka1}
Suppose 
\begin{equation}\label{eqn:ka1}
ka>n-k+\frac{k^2}n.
\end{equation}
Then {\em(\ref{eqn:R})} holds for $1\le q\le k-1$.
\end{lemma}
\begin{proof}
By Lemma \ref{lempg}, it is sufficient to prove that
$$
n-q(a+1)+ \sqrt{(n-q(a+1))^2 + 4q(k-q)}<2n-\frac{2qn}k,
$$
i.~e.
\begin{equation}\label{eqn:sqrt}
 \sqrt{(n-q(a+1))^2 + 4q(k-q)}<n+q(a+1)-\frac{2qn}k.
\end{equation}

We show first that the right-hand side of (\ref{eqn:sqrt}) is positive. In fact
$$
n+q(a+1)-\frac{2qn}k>0\Leftrightarrow k(a+1)>2n-\frac{nk}q.
$$
Now $\frac{nk}q>n$, while $k(a+1)\ge n$ by (\ref{eqn:ka}). So 
$k(a+1)>2n-\frac{nk}q$ as required. 

The inequality (\ref{eqn:sqrt}) is therefore 
equivalent to
$$
(n-q(a+1))^2+4q(k-q)<(n+q(a+1))^2-\frac{4qn}k(n+q(a+1)) +\frac{4q^2n^2}{k^2},
$$
i.~e.
\begin{eqnarray*}
4q(k-q)&<&4nq(a+1)-\frac{4qn}k(n+q(a+1))+\frac{4q^2n^2}{k^2}\\
&=&4q(k-q)\left(\frac{n}k(a+1)-\frac{n^2}{k^2}\right).
\end{eqnarray*}
Dividing by $4q(k-q)$ and rearranging, this becomes (\ref{eqn:ka1}).
\end{proof}
 
In view of Lemma \ref{lemka1}, we now need to deal only with the cases 
$ka=n-k+\frac{k^2-1}n$ and $ka=n-k+\frac{k^2}n$. The second case is 
impossible since $n>k$. It remains to consider the case
\begin{equation}\label{eqn:ka2}
ka=n-k+\frac{k^2-1}n,
\end{equation}
which gives
\begin{equation}\label{eqn:ka3}
n-q(a+1)=n-\frac{qn}k-\frac{q(k^2-1)}{kn}.
\end{equation}
\begin{lemma}\label{lemka2}
Suppose {\em (\ref{eqn:ka3})} holds and $1\le q\le k-1$. Then
$$
(n-q(a+1))^2+4q(k-q)<\left(n-\frac{qn}k+\frac{q(k^2-1)}{kn}+\frac2k\right)^2.
$$
\end{lemma}
\begin{proof}
We need to show that
\begin{equation}\label{ eqn:ka4}
4q(k-q)<\left(n-\frac{qn}k+\frac{q(k^2-1)}{kn}+\frac2k\right)^2-
\left(n-\frac{qn}k-\frac{q(k^2-1)}{kn}\right)^2.
\end{equation}
Now the right-hand side of (\ref{ eqn:ka4}) is equal to
$$
\begin{array}{l}
\displaystyle{\left(2n-\frac{2qn}k+\frac2k\right)\left(\frac{2q(k^2-1)}{kn}+\frac2k\right)}\\
\quad\quad=\displaystyle{\frac4{k^2}\left((k-q)q(k^2-1)+(k-q)n+\frac{q(k^2-1)}n+1\right)}.
\end{array}
$$
So we need to show that
$$
0<-q(k-q)+n(k-q)+\frac{q(k^2-1)}n+1,
$$
which is true since $n>k>q$.
\end{proof}
 
Suppose now that (\ref{eqn:ka2}) holds. Then Lemmas \ref{lempg} 
and \ref{lemka2} imply that
\begin{eqnarray*}
R(n,a,k,q)&\le&\left\lfloor\frac12\left(n-\frac{qn}k-\frac{q(k^2-1)}{kn}+
n-\frac{qn}k+\frac{q(k^2-1)}{kn}+\frac2k\right)\right\rfloor\\
&=&\left\lfloor n-\frac{qn}k+\frac1k\right\rfloor.
\end{eqnarray*}
Note that, if $n-\frac{qn}k+\frac1k$ is an integer, then Lemma \ref{lemka2} 
implies that this inequality is strict. Hence in all cases
$$
R(n,a,k,q)\le n-\frac{qn}k.
$$
Finally $\gcd(n,k)=1$ by (\ref{eqn:ka2}) and $0<q<k$, so $\frac{qn}k$ is 
not an integer. Hence (\ref{eqn:R}) holds. This completes the proof of the theorem.

\end{proof}

\section{The general case}

We now start on the computations of $C_{12}$ and $C_{21}$,
 where we continue to assume that $0<k<n$. With the notation 
of sections 1 and 2, note that, by \cite[Lemma 6.5]{bgn}, the flip loci at 
any critical value can be constructed using only those critical data sets 
for which there exist $(E_1,V_1)$ and $(E_2,V_2)$ which are both $\alpha$-stable either for 
$\alpha=\alpha_c^-$ or for $\alpha=\alpha_c^+$. Since we prefer to have 
purely numerical conditions on our critical data sets, we shall say that 
$A_c=(\alpha_c,n_1,d_1,k_1,n_2,d_2,k_2)$ is {\em allowable} if the numerical
 conditions (\ref{eqn1.3}) and (\ref{eqn:alphac}) hold together with the 
Brill-Noether conditions
\begin{equation}\label{eqn:bn1}
d\ge \frac 1k(n^2-1)-(n-k),\quad 
d_1\ge \frac 1{k_1}(n_1^2-1)-(n_1-k_1)
\end{equation}and
\begin{equation}\label{eqn:bn2}
\mbox{either}\ k_2=0, n_2=1\ \mbox{or}\ 
k_2\ge1, d_2\ge \frac 1{k_2}(n_2^2-1)-(n_2-k_2).
\end{equation}

\begin{proposition} \label{prop2.1}
Let $A_c$ be an allowable critical data set with $k_2=0$. Then $C_{12} > 0$ and $C_{21} >0$.
\end{proposition}

\begin{proof}
By (\ref{eqn:bn2}), we have $n_2=1$ and thus $n_1=n-1$.
Now (\ref{eqn1.4}) and (\ref{eqn4a}) imply
$$
\alpha_c = \frac{ne+t}{k} \quad \mbox{with} \quad 0<e<l.
$$
So by (\ref{eqn1.1}) and (\ref{eqn3a})
\begin{eqnarray*}
C_{12} & =& -(n-1) -ne -t +k(a+e+1)\\
& =& ka -(n-k)e -t -n+ k+1\\
& =& l(n-k) +t +m -(n-k)e -t -n +k+1\\
& =& (l-e-1)(n-k) + m +1  > 0
\end{eqnarray*}
and by (\ref{eqn1.2}) and (\ref{eqn3a})
$$
C_{21} = -(n-1) + ne + t >0.
$$
\end{proof}

\begin{cor}\label{cornon}
If $G(\alpha;n,d,1)$ is non-empty for some $\alpha$ with $t < \alpha < \frac{d}{n-1} - \frac{mn}{n-1}$, then it is non-empty for all such $\alpha$.
\end{cor}

\begin{proof}
For $k=1$, (\ref{eqn1.3}) implies that $k_2=0$ for all critical data sets. 
Hence Proposition \ref{prop2.1} and Corollary \ref{cor:birat} imply the assertion.
\end{proof}

This was proved by a different method in \cite{ln}.

Another case that can be handled easily is when $k_1\ge n_1$.

\begin{proposition} \label{lemma2.3}
$C_{12} >0$ for any allowable critical data set 
with $k_1 \geq n_1$.
\end{proposition}

\begin{proof} By (\ref{eqn1.1})
$$
C_{12} = (k_1-n_1)(n_2+d_2) + d_1n_2 - k_1k_2.
$$
Now $k<n$ implies $k_2<n_2$, so
$$
C_{12} > (k_1-n_1)(n_2+d_2) + (d_1 - k_1)n_2. 
$$
Hence it suffices to show that $d_1 \geq n_1$, since then
$$
C_{12} > (k_1-n_1)(n_2+d_2) + (n_1-k_1)n_2 = d_2(k_1-n_1), 
$$
which is non-negative, since $d_2>0$ by (\ref{eqn1.3}).

In order to see that $d_1 \geq n_1$, suppose first that $k_1 = n_1 + \nu$ 
with $\nu \geq 1$. Then (\ref{eqn:bn1}) implies that
$$
d_1 \geq \frac{1}{n_1+\nu}(n_1^2-1) + \nu \geq n_1.
$$ 
If $n_1=k_1 \geq 2$, the same result gives $d_1 \geq n_1 - \frac{1}{n_1}$
which implies the assertion, since $d_1$ is an integer. Finally, 
if $n_1=k_1=1$ and $d_1<1$, then 
$d_1 = 0$, which implies $\alpha_c = \frac{d}{n-k}$. This contradicts 
(\ref{eqn:alphac}).
\end{proof}

In view of these propositions, we now assume that $k_2\ge1$ and $k_1<n_1$. 
For this case, we need to rearrange the formula for $C_{12}$. We have, 
using (\ref{eqn3a}), (\ref{eqn1.4}) and (\ref{eqn1.5}),
\begin{eqnarray*}
C_{12} &=& -n_1n_2 - (ne + n_2t) +k_1(n_2a+e) + k_1n_2-k_1k_2\\
&= &-n_1n_2 -(n-k_1)\left(-\frac{k_2}{k}t + l\left(n_2- \frac{k_2}{k}n\right) - \frac{f}{k}\right) -n_2t \\
& & +\frac{k_1n_2}{k}(l(n-k)+m+t) +k_1n_2 -k_1k_2.
\end{eqnarray*}
Hence
\begin{eqnarray*}
kC_{12} &=& l[-(n-k_1)(n_2k-nk_2) + k_1n_2(n-k)]\\
&&+ t[(n-k_1)k_2 -n_2k +k_1n_2]\\
&&+ (n-k_1)f + k_1n_2m+ k(k_1n_2 -k_1k_2-n_1n_2)\\
& =& nk_2(n_1-k_1)l + k_2(n_1-k_1)t  + (n-k_1)f\\
&& + k_1n_2m+ k(k_1n_2 -k_1k_2-n_1n_2)
\end{eqnarray*}
and thus 
\begin{equation} \label{eqn2.1}
kC_{12} = (n_1-k_1)(nk_2l + k_2t -kn_2) + (n-k_1)f   + k_1n_2m -kk_1k_2.
\end{equation}

We now use the assumption $k_2\ge1$. The condition (\ref{eqn:bn2}) 
is equivalent by 
(\ref{eqn3a}) and (\ref{eqn1.5}) to
\begin{eqnarray*}
k(n_2a+e) &=& n_2(l(n-k)+m+t) - k_2t + kn_2l - nk_2l -f \\
&\geq& \frac{k}{k_2}(n_2^2-1) -k(n_2-k_2)
\end{eqnarray*}
and thus to
\begin{equation} \label{eqn2.2}
n_2m\ge-(n_2-k_2)(ln+t+k)+f+\frac{k}{k_2}(n_2^2-1).
\end{equation}

We can now prove a partial result for $k_1<n_1$, which will be sufficient 
for our purposes.

\begin{lemma}\label{lemk1n1}
$C_{12}>0$ for any allowable critical data set with $k_1<n_1$, $kk_1<nk_2$ 
(resp. $kk_1\le nk_2$) and $nk_2l+k_2t-kn_2\le0$ (resp. $nk_2l+k_2t-kn_2<0$).
\end{lemma}

\begin{proof}
Inserting (\ref{eqn2.2}) in (\ref{eqn2.1}), we get
\begin{eqnarray*}
kC_{12} &\geq& (n_1-k_1)(nk_2l + k_2t -kn_2) + (n-k_1)f -kk_1k_2\\
& &- (n_2-k_2)k_1(ln + t + k) +k_1f + \frac{kk_1}{k_2}(n_2^2 -1)\\
& = &(k_2(n_1-k_1) - k_1(n_2-k_2))(nl+t) -kn_2(n_1-k_1)\\
&& - kk_1(n_2-k_2)  -kk_1k_2 + \frac{kk_1}{k_2}(n_2^2 -1) + (n-k_1)f +k_1f 
\end{eqnarray*}
i.~e.
\begin{equation} \label{eqn2.3}
kk_2C_{12} \geq   (k_2n_1-k_1n_2)(nk_2l+k_2t-kn_2) + k_2nf - kk_1.
\end{equation}
Note that $k_2n_1-k_1n_2<0$ and $f\ge1$. The result follows.\\
\end{proof}

The formulae (\ref{eqn2.1}) and (\ref{eqn2.3}) are complementary to one 
another in that the first is of use when $nk_2l+k_2t-kn_2\ge0$ and the 
second when $nk_2l+k_2t-kn_2\le0$. This is sufficient to handle another 
special case.

\begin{proposition} \label{prop2.4}
Let $A_c$ be an allowable critical data set with $k_1 =1$. Then $C_{12} > 0$.
\end{proposition}

\begin{proof}
If $k_2=0$, this follows from Proposition \ref{prop2.1}, while, if $n_1=1$, 
it follows from Proposition \ref{lemma2.3}.
If $n_1 \geq 2$, $k_2\ge1$ and $nk_2l+k_2t - kn_2 \leq 0$, then 
we have $kk_1 = k < n \leq nk_2$ and the result follows from Lemma 
\ref{lemk1n1}. 

If $nk_2l+k_2t-kn_2 >0$, then (\ref{eqn2.1}) gives 
$$
kC_{12} \geq (n_1-1) + (n-1)f +n_2m - kk_2.
$$
From (\ref{eqn6}) we get $f \equiv t + nl  \equiv -m \; \mbox{mod}\, k$; moreover $f \geq 1$. 
If $0 \leq m \leq k-1$, then $f \geq k-m$ and 
$$
kC_{12} \geq (n_1-1) + (n-1)(k-m) + n_2m - kk_2 = n_1-1 + (n-k)k - (n_1-1)m. 
$$
But $m<k$ and $n_1-1 < n-1-k_2 = n-k$, since $n_2 > k_2$ by (\ref{eqn1.3}). 
So $kC_{12}>0$. 

Finally, if $m\ge k$, then $kC_{12} \geq n_2k -kk_2 >0$.
\end{proof}

We now turn to look at $C_{21}$.
\begin{lemma}\label{lemc21}
Suppose $k_2\ge1$. Then $C_{21}>0$ in each of the following cases:
\begin{itemize}
\item[(i)] $e\ge1$, $k_2\ge2$, $n\ge k_2(k_1+1)$,
\item[(ii)] $e\ge1$, $k_2=1$, $n\ge2k_1+1$, 
\item[(iii)] $e\le0$, $n\ge k(1+k_1k_2)$.
\end{itemize}\end{lemma}

\begin{proof}
Substituting $d_1=d-d_2$ in (\ref{eqn1.2}) and using (\ref{eqn3a}),
\begin{eqnarray}
C_{21}&=&-n_1n_2+d_2(n-k_2)-d(n_2-k_2)+k_2(n_1-k_1)\nonumber\\
&=&-n_1n_2+(n_2a+e)(n-k_2)-(na-t)(n_2-k_2)+k_2(n_1-k_1)\nonumber\\
&=& n_1(k_2(a+1)-n_2)+e(n-k_2)+t(n_2-k_2)-k_1k_2\label{eqn:c211}
\end{eqnarray}
By (\ref{eqn3a}) and (\ref{eqn:bn2}), we have
$$
k_2(a+1)-n_2\ge\frac{k_2^2-1-k_2e}{n_2}.
$$
So
\begin{eqnarray}
C_{21}&\ge&\frac{n_1(k_2^2-1-k_2e)}{n_2}+e(n-k_2)+t(n_2-k_2)-k_1k_2\nonumber\\
&=&(n_2-k_2)\left(\frac{ne}{n_2}+t\right)+\frac{n_1(k_2^2-1)}{n_2}-k_1k_2\label{eqn:c212}
\end{eqnarray}

If $e\ge1$, this gives
\begin{eqnarray*}
C_{21}&\ge&(n_2-k_2)\frac{n}{n_2} +\frac{n_1(k_2^2-1)}{n_2}-k_1k_2\\
&=& n-k_2(k_1+1)+\frac{n_1(k_2^2-k_2-1)}{n_2}.
\end{eqnarray*}
So $C_{21}>0$ if $n\ge k_2(k_1+1)$ and $k_2\ge2$, proving (i). 
If $k_2=1$, (\ref{eqn1.3}) gives $\frac{n_1}{n_2}<k_1$, so
$$
C_{21}\ge n-(k_1+1)-\frac{n_1}{n_2}>n-2k_1-1.
$$
So $C_{21}>0$ if $n\ge 2k_1+1$, proving (ii). 

If $e\le0$, then $\frac{ne}{n_2}\ge\frac{ke}{k_2}$, so 
$\frac{ne}{n_2}+t\ge\frac{ke}{k_2}+t\ge\frac1{k_2}$ by (\ref{eqn4a}). So 
(\ref{eqn:c212}) gives
\begin{eqnarray}
C_{21}&\ge&\frac{n_2}{k_2}-1+\frac{n_1(k_2^2-1)}{n_2}-k_1k_2\label{eqn:c213}\\
&>&\frac nk+\frac{n_1(k_2^2-1)}{n_2}-(1+k_1k_2).\nonumber
\end{eqnarray}
So $C_{21}>0$ if $n\ge k(1+k_1k_2)$, proving (iii).
\end{proof}

\begin{rem}\label{rem6}
{\em If $k_2=1$, $e\le0$, then (\ref{eqn:c213}) gives $C_{21}\ge n_2-(k_1+1)=n_2-k$, 
with equality possible only if $e=0$, $t=1$.}
\end{rem}

These results are not sufficient for us to determine precisely when 
$C_{12}>0$ or $C_{21}>0$. We shall see in sections 7 and 10 that 
both $C_{12}$ and $C_{21}$ can be $0$. However we can now prove

\begin{theorem}\label{thlge}
Let $k$ be a fixed positive integer. Then there are only finitely many 
allowable critical data sets with $n>k$ for which $C_{12}\le0$ or $C_{21}\le0$.
\end{theorem}
\begin{proof}
By Proposition \ref{prop2.1}, we can suppose that $k_2\ge1$.

Combining Proposition \ref{lemma2.3} with Lemma \ref{lemk1n1}, we see 
that $C_{12}>0$ when $n>\frac{kk_1}{k_2}$, except possibly when
$$
k_1<n_1\quad\mbox{and}\quad nk_2l+k_2t-kn_2>0.
$$
In this case we apply (\ref{eqn2.1}). Since $f\ge1$ and $m\ge0$, we get
$$
kC_{12}>n-k_1-kk_1k_2.
$$
So $C_{12}>0$ if $n\ge k_1+kk_1k_2$. It remains to show that, if we fix $n$ 
as well as $k$, then $C_{12}>0$ for all but finitely many values of $d$.
In view of Proposition \ref{lemma2.3}, we need only 
prove this when $k_1<n_1$. In this case, it follows 
immediately from (\ref{eqn2.1}) that $C_{12}>0$ for all sufficiently 
large values of $l$, say $l\ge A$. But it follows easily from the 
definition of $l$ that this certainly holds if
$$
kd\ge (n-k)((A+1)n-1).
$$

Turning to $C_{21}$, it follows at once from Lemma \ref{lemc21} that, 
for any fixed $k$, $C_{21}>0$ for all sufficiently large $n$. If we fix $n$ 
as well as $k$, and insert $e>-\frac{k_2}{k}t$ in (\ref{eqn:c211}), we obtain
$$
C_{21}>n_1(k_2(a+1)-n_2)+t\left(n_2-k_2-\frac{nk_2}{k}+\frac{k_2^2}k\right)-k_1k_2.
$$
So $C_{21}>0$ for all sufficiently large values of $a$ and hence for all but 
finitely many values of $d$.
\end{proof}
\begin{cor}\label{corlge1}
Let $k$ be a fixed positive integer. Then, for all but finitely many 
pairs $(n,d)$ with $n>k$, one of the following two possibilities holds:
\begin{itemize}
\item $G(\alpha;n,d,k)=\emptyset$ for all $\alpha${\em;}
\item $G(\alpha;n,d,k)\ne\emptyset$ for all $\alpha$ such that 
$$
\frac tk<\alpha<\frac{ln+t}k.
$$
\end{itemize}
\end{cor}
\begin{proof}
This follows from the theorem, Corollary \ref{cor:birat} and (\ref{eqn:alphac}).
\end{proof}

When $t=0$, we have a stronger result.
\begin{cor}\label{corlge2}
Let $k$ be a fixed positive integer. Then, for all but finitely many pairs 
$(n,a)$ such that $n>k$ and {\em(\ref{eqn:ka})} holds, the moduli space
$G(\alpha;n,na,k)\ne\emptyset$ if and only if
$$
0<\alpha<\frac{ln}k.
$$
\end{cor}
\begin{proof}
This follows from Corollary \ref{corlge1} and Theorem \ref{th:0+}.
\end{proof}

We finish this section by showing how we can use these results to 
construct $\alpha$-stable coherent systems for certain values of $t>0$. 

We begin with a lemma
\begin{lemma}\label{lemext1}
Suppose that $t\ge1$, $ka\ge n-k+t$ and $(E_1,W)$ is a coherent system of type 
$(t,t(a-1),k)$. Then
$$
\dim\mbox{\em Ext}^1((E_1,W),(\cO(a),0))\ge n-t.
$$
\end{lemma}
\begin{proof}
By (\ref{eqn1.1}) and \cite[equation (8)]{bgn},
\begin{eqnarray*}
\dim\mbox{Ext}^1((E_1,W),(\cO(a),0))&\ge&-t-at+t(a-1)+k(a+1)\\
&=&-2t+k(a+1)\\
&\ge&-2t+k+n-k+t=n-t.
\end{eqnarray*}
\end{proof}

\begin{proposition}\label{propt>0}
Suppose $k\ge2$, $ka\ge n-k+t$ and that one of the following four conditions 
holds:
\begin{itemize}
\item $t=1$ and $a\ge k$;
\item $t=k-1$ and $a\ge2$;
\item $t=k$ and $a\ge3$;
\item $t>k$, $ka\ge t+\frac{k^2-1}t$ and $C_{12}>0$ 
for all allowable critical data sets for coherent systems of type 
$(t,t(a-1),k)$.
\end{itemize}
Then
$$
G(( t/k)^+;n,d,k)\ne\emptyset.
$$
\end{proposition}
\begin{proof}
We show first that the hypotheses imply that
$$
G( t/k;t,t(a-1),k)\ne\emptyset.
$$

For $t=1$, we require only the condition $h^0(\cO(a-1))\ge k$, which is equivalent 
to $a\ge k$.

For $t=k-1$, the result follows from \cite[proposition 6.4]{ln}.

For $t=k$, it follows from \cite[Proposition 6.3]{ln} that
$$
G(\tilde{\alpha};t,t(a-1),t)\ne\emptyset
$$
for some $\tilde{\alpha}>0$ if and only if $a\ge3$. Taking a general element
$(E,V)$ of this moduli space, we can assume by \cite[Theorem 3.2 and Proposition 3.6]{ln} that 
$E=\cO(a-1)^t$ and that $V$ generically generates $\cO(a-1)^t$. If $(F,W)$ is a 
coherent subsystem of $(E,V)$ which contradicts $\alpha$-stability for some 
$\alpha>0$, then we must have $F=\cO(a-1)^r$, $\dim W=r$ for some $r$, $0<r<t$. 
But then $(F,W)$ contradicts $\alpha$-stability for all $\alpha>0$ and in 
particular for $\alpha=\tilde{\alpha}$. This is a contradiction, establishing that 
$(E,V)$ is $\alpha$-stable for all $\alpha>0$.

Finally, if $t>k$, the hypothesis on the allowable critical data sets implies, by Theorem \ref{th:0+} 
and Corollary \ref{cor:birat}, that $G( t/k;t,t(a-1),k)\ne\emptyset$ provided that
$$
0<\frac tk<\frac{t(a-1)}{t-k}-\frac{m't}{k(t-k)}
$$
for a certain integer $m'$ with $0\le m'<t-k$. This condition is equivalent to
$$
t(t-k)<kt(a-1)-m't,
$$
i.~e. $ka>t+m'$. But $m'<t-k<n-k$, so this follows from the hypothesis $ka\ge n-k+t$.

We now consider extensions
$$
0\ra(\cO(a)^{n-t},0)\ra(E,V)\ra(E_1,W)\ra0,
$$
where $(E_1,W)$ is a $t/k$-stable coherent system of type $(t,t(a-1),k)$. Note that
$$
\mu_{t/k}(E_1,W)=a-1+\frac tk\cdot\frac kt=a.
$$
By Lemmas \ref{lemext} and \ref{lemext1}, the general extension of this form is 
$(t/k)^+$-stable. This completes the proof.
\end{proof}
\begin{cor}\label{cort>k}
Let $k$ be a fixed integer, $k\ge2$. For all but a finite number of pairs $(n,d)$ for which $n>k$, 
$ka\ge n-k+t$ and one of 
the conditions
\begin{itemize}
\item $t=1$ and $a\ge k$,
\item $t=k-1$ and $a\ge2$,
\item $t=k$ and $a\ge3$,
\item $t>k$ and $ka\ge t+\frac{k^2-1}t$
\end{itemize}
holds, the moduli space $G(\alpha;n,d,k)\ne\emptyset$ 
if and only if
$$
\frac tk<\alpha<\frac d{n-k}-\frac{mn}{k(n-k)}.
$$
\end{cor}
\begin{proof}
In view of Theorem \ref{thlge}, we can assume that $C_{12}>0$ for all allowable
critical data sets for coherent systems of type $(n,d,k)$. In the case $t>k$, 
a given pair $(t,a)$ can arise from only finitely many pairs $(n,d)$ which satisfy the 
condition $ka\ge n-k+t$. We can therefore also assume that $C_{12}>0$ for all allowable
critical data sets for coherent systems of type $(t,t(a-1),k)$. The proposition now implies that
$G((t/k)^+;n,d,k)\ne\emptyset$ and the result follows from Corollary \ref{cor:birat}.
\end{proof}

\section{The case $k=2$}

In the case $k=2$, we can use the methods developed above to give a simpler 
proof of \cite[Theorem 5.4]{ln}. Note first that it follows from Propositions 
\ref{prop2.1} and \ref{prop2.4} that $C_{12}>0$ for any allowable critical data set 
and that $C_{21}>0$ except possibly when $k_1=k_2=1$.

\begin{lemma} \label{lemma3.2}
Let $n\ge3$ and let $A_c$ be an allowable critical data set with $k_1=k_2=1$. Then $C_{21} > 0$.
\end{lemma}

\begin{proof}
If $e\ge1$, this follows at once from Lemma \ref{lemc21}(ii).

If $e\le0$, Remark \ref{rem6} gives $C_{21}\ge n_2-2$, with equality possible only if 
$e=0$ and $t=1$. Now $n_2\ge2$ by (\ref{eqn1.3}). Hence $C_{21}>0$ except 
possibly when $e=0$, $t=1$, $n_2=2$, and then $n_1=1$ by (\ref{eqn1.3}). 
Moreover $d=3a-1$ and (\ref{eqn:bn1}) implies that $d \geq 3$. Hence $a \geq 2$ and 
by (\ref{eqn1.2}) and (\ref{eqn3a})
$$
C_{21} = -2 + 2a-2(a-1)+(a-1+1-1)=a-1 \geq 1,
$$
which completes the proof of the lemma.
\end{proof}

\begin{cor} \label{cor3.3}
If $G(\alpha;n,d,2)$ is non-empty for some $\alpha$ with 
$\frac{t}{2} < \alpha < \frac{d}{n-2} - \frac{mn}{2(n-2)}$, 
then it is non-empty for all such $\alpha$. 
\end{cor}

\begin{proof}
It suffices to show that $C_{12}$ and $C_{21}$ are both positive for all
allowable critical data sets $A_c$. But for $k=2$ either $k_2=0$
or $k_1=k_2=1$. Hence Propositions \ref{prop2.1} and \ref{prop2.4} 
and Lemma \ref{lemma3.2} imply the assertion.
\end{proof}

For the proof of the full result of \cite{ln} for $k=2$, it remains to 
determine when there exists an $\alpha$-stable coherent system
of type $(n,d,2)$ for some $\alpha$. 

\begin{proposition} \label{prop3.6}
Suppose $n \geq 3$,  $l \geq 1$, $d \geq \frac{1}{2}n(n-2) + \frac{3}{2}$ 
and  $(n,d) \neq (4,6)$. Then there exists a $(t/2)^+$-stable coherent 
system $(E,V)$ of type $(n,d,2)$.
\end{proposition} 

\begin{proof}
For $t=0$, this has already been proved in Theorem \ref{th:0+}.

For $t \geq 1$, it is sufficient to verify that the conditions of Proposition
 \ref{propt>0} are satisfied. Note that the 
hypothesis $l\ge1$ is equivalent to
\begin{equation}\label{eqn:2a1}
2a\ge n-2+t,
\end{equation}
while the Brill-Noether condition $d\ge\frac12(n^2-1)-(n-2)$ is easily 
seen to be equivalent to
\begin{equation}\label{eqn:2a2}
2a\ge n-2+\frac{3+2t}n.
\end{equation}

For $t=1$, (\ref{eqn:2a2}) gives $2a\ge n-2+\frac5n$, which implies $a\ge2$ as 
required. For $t\ge3$, the condition $2a\ge t+\frac3t$ follows 
from (\ref{eqn:2a1}), while $C_{12}>0$ holds always for $k=2$. 
For $t=2$, (\ref{eqn:2a1}) gives $2a\ge n$, which 
implies $a\ge3$ as required except for $n=3,4$. We are left therefore 
with the two cases $(n,d)=(3,4)$ and $(n,d)=(4,6)$, for both of which 
$a=2$. The case $(n,d)=(4,6)$ has been excluded in the statement, 
so we need only to prove the proposition for $(n,d)=(3,4)$.

In this case, we have $a=t=2$. The moduli space $G(1;2,2,2)$ is empty by \cite[Proposition 5.6]{ln}, 
but there do exist $1$-semistable coherent systems of type $(2,2,2)$, 
which have the form
$$
(E_1,W)=(\cO(1),W_1)\oplus(\cO(1),W_2).
$$
Since $h^0(\cO(1))=2$, we can take $W_1$ and $W_2$ to be distinct 
subspaces of $H^0(\cO(1))$ of dimension $1$. Let $(E,V)$ be the 
general extension of the form
$$
0 \ra (\cO(2),0) \ra (E,V) \ra (E_1,W) \ra 0.
$$
Comparing this with (\ref{eqn:nots}), it is easy to verify (\ref{eqn:hyp2}). 
It follows from Lemma \ref{lemext} that $(E,V)$ is $1^+$-stable as required.
\end{proof}

We can now restate \cite[Theorem 5.4]{ln}.
\begin{theorem}\label{th:k=2}
Suppose $n\ge3$. Then $G(\alpha;n,d,2)\ne\emptyset$ for some 
$\alpha$ if and only if $l\ge1$, $d\ge\frac12n(n-2)+\frac32$ and $(n,d)\ne(4,6)$. 
Moreover, when these conditions hold, $G(\alpha;n,d,k)\ne\emptyset$ if and 
only if
$$
\frac t2<\alpha<\frac d{n-2}-\frac{mn}{2(n-2)}.
$$
\end{theorem}
\begin{proof}
The stated conditions are sufficient by Proposition \ref{prop3.6}. Conversely,
if $G(\alpha;n,d,k)\ne\emptyset$, then $l\ge1$ and $d\ge\frac12n(n-2)+\frac32$ 
by \cite[Remark 4.3 and Corollary 3.3]{ln}. It is easy to prove that there do 
not exist $\alpha$-stable coherent systems of type $(4,6,2)$ (see the first 
paragraph of the proof of \cite[Theorem 5.4]{ln}). For the last part, see 
Corollary \ref{cor3.3}.
\end{proof}

\section{The case $k=3$}

Now suppose $k=3$. In this section we will show that $C_{12}$ is 
positive for all allowable critical data sets $A_c$ and determine those 
$A_c$ for which $C_{21} = 0$. As a consequence, we give examples for which 
the lower bound of \cite[Proposition 4.1]{ln} for those $\alpha$, for which there exist $\alpha$-stable 
systems, is not best possible.

\begin{proposition} \label{prop4.1}
Let $n\ge4$. Suppose $A_c$ is an allowable critical data set with $k=3$. Then $C_{12} >0$.
\end{proposition} 

\begin{proof}
The cases $k_1=3$ and $k_1 = 1$ are covered by Propositions 
\ref{prop2.1} and \ref{prop2.4}. So suppose $k_1=2,\; k_2=1$.
Then (\ref{eqn1.3}) implies
$$
n_1< 2n_2,
$$
and by (\ref{eqn6}) 
$$
f \equiv m \; \mbox{mod} \;3.
$$
According to (\ref{eqn2.1}) 
$$
3C_{12} = (n_1-2)(nl+t-3n_2) +(n-2)f + 2n_2m -6.
$$
If $n_1 \leq 2$, $C_{12}$ is positive by Proposition \ref{lemma2.3}. 
So let $n_1 \geq 3$. If $nl+t-3n_2 \geq 0$, then $C_{12} >0$, since
$n_2 \geq 2$ and thus $n \geq 5$ and either $f$ and $m$ are both positive 
or $f\ge3$. 

If $nl+t-3n_2 < 0$, then, using (\ref{eqn2.3}), we get
$$
3C_{12} \geq (2n_2-n_1)(3n_2-nl-t) +nf - 6,
$$
which is positive for $n \geq 6$. For $n=5$, we have $n_1=3, \; n_2=2$ and
$$
3C_{12}= 5l + t +3f +4m - 12.
$$
This is $\geq 0$ since $l>0$ and either $f$ and $m$ are both positive or 
$f\ge3$. Equality occurs if and only if $l=f=m=1$, $t=0$. But then 
(\ref{eqn1.5}) gives $e= 0$, which contradicts (\ref{eqn4a}).
\end{proof}

\begin{lemma} \label{lem4.3}
Let $n\ge4$. Suppose $A_c$ is an allowable critical data set 
with $k_1 = 1, k_2 = 2$. Then $C_{21} >0$.
\end{lemma}

\begin{proof}
For $e\ge1$, this follows at once from Lemma \ref{lemc21}(i).

For $e\le0$, we need to analyse (\ref{eqn:c212}) and (\ref{eqn:c213}) 
more carefully. In our case (\ref{eqn:c213}) becomes
$$
C_{21}\ge\frac{n_2}2+3\frac{n_1}{n_2}-3.
$$
This is positive if $n_2\ge5$. For $n_2=4$, we note that (\ref{eqn1.3}) implies 
that $n=5$, $n_1=1$. By (\ref{eqn4a}), we have $e\ge-\frac23t+\frac13$, so by (\ref{eqn:c212})
$$
C_{21}\ge2\left(\frac{5e}4+t\right)+\frac34-2
=\frac{5e}2+2t -\frac54
\ge\frac t3+\frac56-\frac54.
$$
This is positive if $t\ge2$. Since $e\le0$, the only remaining case is when 
$t=1$, which implies $e\ge0$, so $C_{21}\ge2t-\frac54>0$.

It remains to consider the case $e\le0$, $n_2=3$. In this case, (\ref{eqn1.3}) gives 
$n=4$, $n_1=1$. We now have by (\ref{eqn:c211})
$$
C_{21}= 2(a+1)-3+2e+t-2=-3+2a+t+2e.
$$
Since $e\le0$, we have $t\ge1$. Moreover (\ref{eqn:bn1}) gives $d\ge4$, so $a\ge2$. 
If $t=1$, then $e=0$, while, if $t=2$ or $t=3$, then $e\ge-1$. So in all cases $C_{21}>0$.
\end{proof}

\begin{proposition} \label{lem4.4}
Suppose $n \geq 4$ and $A_c$ is an allowable critical data set with $k_1 = 2, k_2 = 1$. Then
$C_{21} > 0$, except in the cases 
\begin{itemize}
\item[(a):] $(n_1,n_2,d_1,d_2) = (4,3,7,6), \quad \alpha_c = \frac{3}{2},$

\item[(b):] $(n_1,n_2,d_1,d_2) = (3,3,5,6), \quad \alpha_c = 1,$

\item[(c):] $(n_1,n_2,d_1,d_2) = (2,3,3,6), \quad \alpha_c = \frac{3}{4},$

\item[(d):] $(n_1,n_2,d_1,d_2) = (1,3,1,6), \quad \alpha_c = \frac{3}{5},$
\end{itemize}
where $C_{21}=0$.
\end{proposition}

\begin{proof}
In this case (\ref{eqn1.3}) gives 
\begin{equation} \label{eqn16}
n < 3n_2;
\end{equation}
since $n\ge4$, this implies that $n_2\ge2$.
By (\ref{eqn:c212}) we have
\begin{equation} \label{eqn17}
C_{21} \geq (n_2-1)(\frac{ne}{n_2} + t) -2.
\end{equation}
We distinguish several cases:\\

{\it Case 1}: $n_2 = 2$. According to (\ref{eqn16}), $n_1 = 2$ or 3. Suppose first
$(n_1,n_2)=(2,2)$. By (\ref{eqn1.2}) we have $C_{21} = 2d_2 - d_1 -4$.
By (\ref{eqn:bn1}), $d_1 \geq 2$ and (\ref{eqn1.3}) implies that 
$d_2 \geq d_1+1$. Hence, if $d_1 \geq 3$ or $d_1=2, d_2 \geq 4$, 
we have $C_{21} >0$. In 
the remaining case $(d_1,d_2) = (2,3)$, we have $t =3$, $e=-1$, so 
(\ref{eqn4a}) fails.

Now suppose $(n_1,n_2) =(3,2)$. Then $C_{21} = 3d_2 - d_1 -5$.
By (\ref{eqn:bn1}), $d_1 \geq 3$ and (\ref{eqn1.3}) implies $d_2 > \frac{2}{3}d_1$.
Hence, if $d_1 \geq 5$ or $d_1 = 3, d_2 \geq 3$ or $d_1 = 4, d_2 \geq 4$, we have $C_{21} > 0$. 
In the remaining case $(d_1,d_2) = (4, 3)$, we have $t=3$, $e=-1$, 
so again (\ref{eqn4a}) fails.\\

{\it Case 2}: $e \geq 1$, $n_2 \geq 3$. (\ref{eqn17}) gives
$C_{21} \geq (n_2-1)(1+\frac{n_1}{n_2}) -2 >0$.\\

{\it Case 3}: $e\le0$ $n_2 \geq 3$, $(n_2,t,e)\ne(3,1,0)$. The result follows from Remark 
\ref{rem6}.\\

 {\it Case 4}: $e = 0$, $t=1$, $n_2 = 3$. 
By (\ref{eqn16}) we have
$1 \leq n_1 \leq 5$. Moreover $d_2=3a$ by (\ref{eqn3a}) and hence 
$d_1=n_1a-t=n_1a-1$. So, by (\ref{eqn1.2}),
$$
C_{21}=n_1d_2-2d_1-2(n_1+1)=n_1(a-2).
$$
By (\ref{eqn:bn2}), $d_2\ge6$, so $a\ge2$ and $C_{21}\ge0$. 
Now the Brill-Noether inequality $d\ge\frac13(n^2-1)-(n-3)$ gives
\begin{equation}\label{eqn:3a}
3a\ge n-3+\frac{8+3t}n=n-3+\frac{11}n.
\end{equation}
Using this, we see that $a=2$ only in the four cases listed 
(note that $n_1=5$ does not occur, since then (\ref{eqn:3a}) gives $a\ge3$). 
One can easily compute $\alpha_c$ in each of the exceptional cases and 
check that (\ref{eqn:alphac}) holds.
\end{proof}

\begin{proposition} \label{prop4.7}
For all cases other than those covered by Proposition {\em\ref{lem4.4} (a)--(d)}, 
if $G(\alpha;n,d,3)$ is non-empty for some $\alpha$ with $\frac{t}{3} < \alpha < \frac{d}{n-3} - \frac{mn}{3(n-3)}$,
then it is non-empty for all such $\alpha$. 
\end{proposition}

\begin{proof}
This follows from Propositions \ref{prop2.1} and \ref{prop4.1}, 
Lemma \ref{lem4.3} and Proposition \ref{lem4.4}, together with Corollary \ref{cor:birat}. 
\end{proof}

\section{Existence for $k=3$}

We consider first the existence of $\alpha$-stable coherent systems in the 
exceptional cases of Proposition \ref{lem4.4}.

\begin{proposition}\label{propexc}
In each of the following cases, we have $G(\alpha;n,d,3)=\emptyset$ for 
$\alpha\le\alpha_c$ and $G(\alpha_c^+;n,d,3)\ne\emptyset$:
\begin{itemize}
\item[(a):] $(n,d) = (7,13), \quad \alpha_c = \frac{3}{2},$

\item[(b):] $(n,d) = (6,11), \quad \alpha_c = 1,$

\item[(c):] $(n,d) = (5,9), \quad \alpha_c = \frac{3}{4},$

\item[(d):] $(n,d) = (4,7), \quad \alpha_c = \frac{3}{5}.$
\end{itemize}
\end{proposition}
\begin{proof}
In each case $C_{12}>0$ for all allowable critical data sets by 
Proposition \ref{prop4.1} and, by Lemma \ref{lem4.3} and 
Proposition \ref{lem4.4}, $C_{21}>0$ except for a unique critical 
data set as given in Proposition \ref{lem4.4}. In view of Corollaries 
\ref{cor:birat} and \ref{corflip}, and Remarks \ref{rem2} and \ref{rem3}, 
it is therefore sufficient to prove that there exist $\alpha_c$-stable 
coherent systems $(E_1,V_1)$ of type $(n_1,d_1,2)$ and $(E_2,V_2)$ of 
type $(n_2,d_2,1)$, where $n_1$, $n_2$, $d_1$, $d_2$ are as given in 
Proposition \ref{lem4.4}.

For $(E_2,V_2)$, we have $(n_2,d_2)=(3,6)$ in every case and, with the 
obvious notation, $t_2=0$, $m_2=0$. So, by \cite[Theorem 5.1]{ln}, 
we require
$$
0<\alpha_c<\frac62=3,
$$
which is true in every case.

For $(E_1,V_1)$, in cases (a) and (b) we need to apply Theorem \ref{th:k=2} (or 
\cite[Theorem 5.4]{ln}). Certainly $(n_1,d_1)\ne(4,6)$ and it is easy to check 
that $l_1\ge1$ and $d_1\ge\frac12n_1(n_1-2)+\frac32$; in fact the latter 
was one of the conditions for an allowable critical data set. It remains to prove 
that in each case
$$
\frac{t_1}2<\alpha_c<\frac{d_1}{n_1-2}-\frac{m_1n_1}{2(n_1-2)}
$$
In fact the numbers in each case are given by
\begin{itemize}
\item[(a):] $n_1=4$, $d_1=7$, $t_1=1$, $m_1=1$,
\item[(b):] $n_1=3$, $d_1=5$, $t_1=1$, $m_1=0$,
\end{itemize}
and the result is clear.

In case (c), we have $n_1=k_1=2$, $d_1=3$, so the result follows 
from \cite[Proposition 5.6]{ln}. Finally, in case (d), we have $n_1=1$, $k_1=2$, $d_1=1$, 
so $(E_1,V_1)\simeq(\cO(1),H^0(\cO(1)))$ is $\alpha$-stable for all $\alpha>0$.
\end{proof}

We turn now to the general case.

\begin{proposition}\label{propk=3}
Suppose $n\ge4$, $l\ge1$, $d\ge\frac13n(n-3)+\frac83$ and
$$
(n,d)\ne(7,13),\ (6,11),\ (6,9),\ (5,9),\ (4,7).
$$
Then there exists a $(t/3)^+$-stable coherent system of type $(n,d,3)$.
\end{proposition}
\begin{proof}
For $t=0$, this has already been proved in Theorem \ref{th:0+}. For $t\ge1$, it is sufficient to 
verify that the conditions of Proposition \ref{propt>0} are satisfied.

Note first that the condition
\begin{equation}\label{eqn:3a1}
3a\ge n-3+t
\end{equation}
is equivalent to $l\ge1$.

For $t=1$, we require $a\ge3$. By (\ref{eqn:3a}), the only cases for 
which $a<3$ are when $a=2$ and $4\le n\le7$, giving rise precisely 
to the exceptional cases $(7,13)$, $(6,11)$, $(5,9)$ and $(4,7)$.

For $t=2$, we require $a\ge2$, which follows at once from (\ref{eqn:3a}).

For $t=3$, we require again $a\ge3$. By (\ref{eqn:3a}), we have
$$
3a\ge n-3+\frac{17}n,
$$
which implies $a\ge3$ except in the cases
$$
(n,d)=(6,9),\ (5,7),\ (4,5).
$$
For $(n,d)=(5,7)$ and $(n,d)=(4,5)$, we consider extensions
$$
0\ra(\cO(2)^{n-3},0)\ra(E,V)\ra\bigoplus_{i=1}^3(\cO(1),W_i)\ra0,
$$
where the $W_i$ are distinct subspaces of $H^0(\cO(1))$ of dimension $1$. We need to check 
the inequalities (\ref{eqn:hyp2}), which in this case give $3>n-3$ and 
$3\ge2$. These are valid, so Lemma \ref{lemext} establishes that the 
general $(E,V)$ is $\alpha_c^+$-stable, where here $\alpha_c=1=\frac t3$.

Finally, for $t\ge4$, we certainly have $C_{12}>0$ for all allowable 
critical data sets for coherent systems of type $(t,t(a-1),3)$ by 
Proposition \ref{prop4.1}. The condition $3a\ge t+\frac8t$ follows 
from (\ref{eqn:3a1}) since we now have $n\ge5$.

This completes the proof.
\end{proof}

\begin{rem}\label{rem5}
{\em  The construction fails for $(n,d)=(6,9)$ because we no longer have 
$3>n-3$. In fact, it is easy to see that $G(1^+;6,9,3)=\emptyset$. Indeed, 
if this is not so, then a general element of $G(1^+;6,9,3)$ has 
$E\simeq\cO(2)^3\oplus\cO(1)^3$. Now $E$ has a unique subbundle $F\simeq\cO(2)^3$. 
If $V_1=V\cap H^0(F)\ne0$, then $(F,V_1)$ is a coherent subsystem of $(E,V)$ 
which contradicts $1^+$-stability. So there is an exact sequence
\begin{equation}\label{eqn:693}
0\ra(\cO(2)^3,0)\ra(E,V)\ra(\cO(1)^3,W)\ra0
\end{equation}
with $\dim W=3$. The homomorphism $W\otimes\cO\ra\cO(1)^3$ is 
not an isomorphism, so there exists a section of $\cO(1)^3$ contained in $W$ 
and possessing a zero. This defines a coherent subsystem $(\cO(1),W_1)$ of 
$(\cO(1)^3,W)$ with $\dim W_1=1$. Now consider the pullback
\begin{equation}\label{eqn:231}
0\ra(\cO(2)^3,0)\ra(E_1,V_1)\ra(\cO(1),W_1)\ra0
\end{equation}
of (\ref{eqn:693}). Such extensions are classified by triples $(e_1,e_2,e_3)$ 
with $e_i\in\mbox{Ext}^1((\cO(1),W_1),(\cO(2),0))$. Note that
$$
\mbox{Hom}((\cO(1),W_1),(\cO(2),0))=0.
$$
Hence, from (\ref{eqn:c12}) and (\ref{eqn1.1}), we see that
$$
\dim\mbox{Ext}^1((\cO(1),W_1),(\cO(2),0))=1.
$$
It follows that, using an automorphism of $\cO(2)^3$, we can suppose that $e_2=e_3=0$. 
This means that (\ref{eqn:231}) is induced from an exact sequence
$$
0\ra(\cO(2),0)\ra(E_2,V_2)\ra(\cO(1),W_1)\ra0.
$$
But now $(E_2,V_2)$ is a coherent subsystem of $(E,V)$ which contradicts 
$\alpha$-stability of $(E,V)$ for all $\alpha$. Hence $G(1^+;6,9,3)=\emptyset$ as asserted. 
It follows from Proposition \ref{prop4.7} that $G(\alpha;6,9,3)=\emptyset$ for all $\alpha$.}
\end{rem}

\begin{theorem} \label{thm6.5}
Suppose $n \geq 4$. Then
$G(\alpha;n,d,3) \neq \emptyset$ for some $\alpha > 0$ if and only if $l \geq 1$, 
$d \geq \frac{1}{3}n(n-3) + \frac{8}{3}$ 
and $(n,d) \neq (6,9)$.
Moreover, when these conditions hold, $G(\alpha;n,d,3)\ne\emptyset$ if and only if
$$
\frac{t}{3} < \alpha < \frac{d}{n-3} - \frac{mn}{3(n-3)},
$$
except for the following pairs $(n,d)$, where the range of $\alpha$ is as stated :
$$
\begin{array}{cccccccc}
 for& (4,7):&\frac{3}{5} < \alpha < 7;&&&for & (5,9):& \frac{3}{4} < \alpha < \frac{11}{3};\\
 for& (6,11):& 1 < \alpha < \frac{7}{3};&&&for & (7,13):& \frac{3}{2} < \alpha < \frac{8}{3}.\\
\end{array}
$$
\end{theorem}
\begin{proof}
The necessity of the conditions follows from \cite[Corollary 3.3 and Remark 4.3]{ln} 
and Remark \ref{rem5}. Sufficiency has been proved in Propositions 
\ref{propexc} and \ref{propk=3}. The assertion about the range of $\alpha$ follows from 
Proposition \ref{prop4.7} except for the exceptional cases, when it is a consequence of 
Propositions \ref{prop4.1} and \ref{propexc} and Corollary \ref{cor:birat}.
\end{proof}

\section{The case $k=3$, $n\le3$}

In this section, we will complete the results for $k=3$ by considering the case $n\le3$. 
It is interesting to note that further exceptional cases arise. We begin with a general 
result, which completes \cite[Proposition 6.3]{ln} in the case $t=0$.

\begin{proposition}\label{propna}
For any $n\ge2$, $G(\alpha;n,na,n)\ne \emptyset$ if and only if $a\ge2$ and $\alpha>0$.
\end{proposition}
\begin{proof}
By \cite[Proposition 6.3]{ln}, $G(\alpha; n,na,n)\ne\emptyset$ for some $\alpha>0$ if and 
only if $a\ge2$ and there is then no upper bound on $\alpha$. In view of 
\cite[Corollary 3.4]{ln}, it is therefore sufficient to prove that $G(0^+;n,na,n)\ne\emptyset$ 
if $a\ge2$. But this is exactly what is shown in the proof of the case $t=k$ of 
Proposition \ref{propt>0}.
\end{proof}

\begin{theorem}\label{thm:smalln}
\ 
\begin{itemize}
\item[(i):] $G(\alpha;1,d,3)\ne\emptyset$ if and only if $d\ge2$ and $\alpha>0$.
\item[(ii):] $G(\alpha;2,d,3)\ne\emptyset$ for some $\alpha$ if and only if $d\ge2$. 
Moreover, if $d\ge2$, $G(\alpha;2,d,3)\ne \emptyset$ for all $\alpha>\frac t3$ except 
in the case $d=3$, when $G(\alpha;2,3,3)\ne\emptyset$ if and only if $\alpha>1$.
\item[(iii):] $G(\alpha;3,d,3)\ne\emptyset$ for some $\alpha$ if and only if $d\ge4$. 
Moreover, if $d\ge4$, $G(\alpha;3,d,3)\ne \emptyset$ for all $\alpha>\frac t3$ except 
in the case $d=5$, when $G(\alpha;3,5,3)\ne\emptyset$ if and only if $\alpha>\frac23$.
\end{itemize}
\end{theorem}
\begin{proof}
(i) follows immediately from the fact that $h^0(\cO(d))\ge3$ if and only if $d\ge2$.

(ii): By \cite[Proposition 6.4]{ln}, $G(\alpha;2,d,3)\ne \emptyset$ for some $\alpha$ 
if and only if $d\ge2$ and there is then no upper bound on $\alpha$. Moreover, 
if $d\ge2$ and $t=0$, then, by the 
same proposition, $G(\alpha;2,d,3)\ne \emptyset$ for all $\alpha>0$. 
This completes the proof when $t=0$.

If $t=1$, it follows from Proposition \ref{propt>0} that $G((1/3)^+;2,d,3)\ne\emptyset$ 
if $a\ge3$, i.~e. $d\ge5$; the result now follows from \cite[Corollary 3.4]{ln}. 
Now suppose $d=3$ and consider a coherent system $(E,V)=(\cO(2)\oplus\cO(1),V)$ 
of type $(2,3,3)$ 
such that $V$ generates $E$; in fact we can take $V$ to be a general subspace of 
$H^0(E)$ of dimension $3$. If $(F,W)$ is any coherent subsystem of $(E,V)$ 
with $\mbox{rk}F=1$, then $\dim W\le1$ (otherwise $E/F$ would not be generated by $V/W$).
 Moreover $\deg F\le2$, so $(E,V)$ is $\alpha$-stable provided
$$
2+\alpha<\frac32+\frac{3\alpha}2,
$$
i.~e. $\alpha>1$. Conversely, if $(E,V)$ is any $\alpha$-stable coherent system of 
type $(2,3,3)$, then $E\simeq\cO(2)\oplus\cO(1)$. Since $h^0(\cO(1))=2$, $(E,V)$ has 
a coherent subsystem of type $(1,2,1)$, which contradicts $\alpha$-stability for 
$\alpha\le1$.

(iii): By \cite[Proposition 6.3]{ln}, $G(\alpha;3,d,3)\ne \emptyset$ for some $\alpha$ 
if and only if $d\ge4$ and there is then no upper bound on $\alpha$. So suppose that $d\ge4$.
If $t=0$, Proposition \ref{propna} gives the result. If $t=1$, $d\ge8$ or $t=2$, 
Proposition \ref{propt>0} implies that 
$G((t/3)^+;3,d,3)\ne\emptyset$, and the result follows from \cite[Corollary 3.4]{ln}.

There remains the case $d=5$. We consider a coherent system 
$(E,V)=(\cO(2)^2\oplus\cO(1),V)$ of type $(3,5,3)$, where we choose $V$ 
so that $V$ generically generates $E$, $\dim V\cap H^0(\cO(2)^2)=1$ and the line subbundle 
generated by a non-zero element of $ V\cap H^0(\cO(2)^2)$ has degree $\le1$. If now $(F,W)$ 
is a coherent subsystem of $(E,V)$ with $\mbox{rk}F=1$, we have either $\dim W=0$, $\deg F=2$ 
or $\dim W\le1$, $\deg F\le1$. In the first case, the $\alpha$-stability condition holds for 
$\alpha>\frac13$, in the second case for all $\alpha>0$. Now suppose $(F,W)$ is a 
coherent subsystem of $(E,V)$ with $\mbox{rk}F=2$. Then either $\dim W\le1$, $\deg F=4$ or
$\dim W\le2$, $\deg F\le3$. In the first case, the $\alpha$-stability condition holds 
for $\alpha>\frac23$, in the second case for all $\alpha>0$. Thus we have shown that 
$(E,V)$ is $\alpha$-stable for $\alpha>\frac23$.

Conversely suppose $(E,V)$ is an $\alpha$-stable coherent system of type $(3,5,3)$. 
Then $E\simeq\cO(2)^2\oplus\cO(1)$ or $\cO(3)\oplus\cO(1)^2$. Since $h^0(\cO(1))=2$, $(E,V)$ has a 
coherent subsystem of type $(2,4,1)$, which contradicts $\alpha$-stability for 
$\alpha\le\frac23$.
\end{proof}

\section{An example for $k=4$}

In this section we give an example of an allowable critical data set 
$A_c$ with $0<k<n$ and $C_{12}=0$. According to our earlier results, we must have 
$k\ge4$ and, by Proposition \ref{lemma2.3}, $k_1<n_1$. Now (\ref{eqn1.3}) 
implies that  $k_2 <n_2$, so $n\ge6$. The minimal possible example therefore 
has $n=6$, $k=4$ and one can check that then $k_1=3$, $k_2=1$, $n_1=4$, 
$n_2=2$. The formula (\ref{eqn2.1}) gives
$$
kC_{12} = (n_1-k_1)(nk_2l + k_2t -kn_2) + (n-k_1)f   + k_1n_2m -kk_1k_2,
$$
i.~e.
$$
4C_{12}=6l+t+3f+6m-20.
$$
Now $l\ge1$ and, by (\ref{eqn6}), $f\equiv m\bmod4$. Since $f\ge1$, either 
$f$ and $m$ are both positive or $f\ge4$. It is now easy to check that the only cases giving $C_{12}=0$ are
$$
l=1, f=m=1, t=5;\quad l=1,f=4,m=0,t=2.
$$
In both cases one can check from (\ref{eqn1.5}) that $e=-1$, which implies by 
(\ref{eqn4a}) that $t\ge5$. Thus we are left with just one case in which a simple 
computation gives $d=7$. Note in this case that the necessary condition for 
$\alpha$-stability from \cite[Propositions 4.1 and 4.2]{ln} is
$$
\frac54<\alpha<\frac{11}4.
$$
Now by (\ref{eqn1.4})
 
$$
\alpha_c=\frac{ne+n_2t}{n_2k-nk_2}=\frac{-6+10}{8-6}=2,
$$
which does lie within the given range. The critical data set itself is 
given by
$$
A_c=(\alpha_c,n_1,d_1,k_1,n_2,d_2,k_2)=(2,4,4,3,2,3,1).
$$
One can check (\ref{eqn:bn1}) and (\ref{eqn:bn2}) to show that $A_c$ is allowable.

\begin{proposition}\label{prop9}

\begin{itemize}
\item[(a):] $G(\alpha;6,7,4)\ne\emptyset\quad\mbox{for} \quad\frac54<\alpha<2$.
\item[(b):] $G(\alpha;6,7,4)=\emptyset\quad\mbox{for} \quad\alpha\ge2$.
\end{itemize}
\end{proposition}
\begin{proof}
We show first that $G(2^+;6,7,4)=\emptyset$ and $G(2^-;6,7,4)\ne\emptyset$. 
First note that
$$
C_{21}=-n_1n_2+d_2n_1-d_1n_2+k_2(d_1+n_1-k_1)=-8+12-8+4+4-3=1.
$$
The result will therefore follow from Corollary \ref{corflip} and Remarks 
\ref{rem2} and \ref{rem3} if we prove the existence of $2$-stable coherent 
systems of types $(4,4,3)$ and $(2,3,1)$. In the first case, we have to 
check the conditions of Theorem \ref{thm6.5}; in the second case, 
those of \cite[Theorem 5.1]{ln}. Both computations are easy.

By Corollary \ref{corflip}, we have $G(2^-;6,7,4)=G_2^-$, so $G(2;6,7,4)=\emptyset$. 
It follows from \cite[Corollary 3.4]{ln} that $G(\alpha;6,7,4)=\emptyset$ if $\alpha\ge2$, 
thus proving (b). For (a), we can apply Proposition \ref{propt>0} to show that 
$G((5/4)^+;6,7,4)\ne\emptyset$. In particular, we must show that $C_{12}>0$ for all 
allowable critical data sets for coherent systems of type $(5,5,4)$. From our argument 
above, we have shown that $C_{12}>0$ always if $k=4$ and $n=5$, so this is clear. The 
result now follows from \cite[Corollary 3.4]{ln}.
\end{proof}


\begin{thebibliography}{CAV}

\bibitem{bgn} S. B. Bradlow, O. Garc\'{\i}a-Prada, V. Mu\~noz and P. E. Newstead: 
\emph{Coherent systems and Brill-Noether theory}. 
Internat. J. Math. 14 (2003), 683-733.


\bibitem{ln} H. Lange and P. E. Newstead: 
\emph{Coherent systems of genus 0}.
Internat. J. Math. 15 (2004), 409-424.



\end{thebibliography}
\end{document}